\documentclass{amsart}

\usepackage{amssymb}

\newcommand{\sd}{\rtimes}
\newcommand{\an}{\mbox{{\rm {\tiny an}}}}
\newcommand{\Aut}{\mbox{{\rm Aut}}}

\newtheorem*{theorem}{Theorem}
\newtheorem*{lemma}{Lemma}
\newtheorem{proposition}{Proposition}

\theoremstyle{remark}
\newtheorem{remark}{Remark}

\begin{document}

\title[Mumford curves with maximal automorphism group]{Mumford curves with maximal \\ automorphism group}
\author{Gunther Cornelissen}
\address{Ghent University \\ Dept.\ of Pure Mathematics and Computer Algebra \\ Galglaan 2 \\
9000 Ghent \\ Belgium}
\curraddr{Max-Planck-Institut f{\"u}r Mathematik \\ Vivatsgasse 7 \\ 53111 Bonn \\ Germany}
\email{gc@cage.rug.ac.be}
\thanks{The first author is Post-doctoral fellow of the the Fund for Scientific Research -- Flanders (FWO-Vlaanderen). This work
was done when he was visiting Kyoto University. The main result of this paper
answers positively to a question posed by T.\ Sekiguchi during the 2000 Kinosaki Symposium 
on Algebraic Geometry.}

\author{Fumiharu Kato}
\address{Kyoto University \\ Faculty of Science \\ Dept.\ of Mathematics \\ Kyoto 606-8502 \\ Japan}
\email{kato@kusm.kyoto-u.ac.jp}


\begin{abstract} It is known that a Mumford curve of genus $g \notin \{5,6,7,8 \}$ over a non-archimedean valued 
field of characteristic $p>0$ has at most $2 \sqrt{g} (\sqrt{g}+1)^2$ automorphisms. In this note,
the unique family of curves which attains this bound, and their 
automorphism group are determined.
\end{abstract}

\maketitle

\section{Introduction}

It is well-known (cf.\ Conder \cite{Conder:90})
that if a compact Riemann surface of genus $g \geq 2$ attains
Hurwitz's bound $84(g-1)$ on its number of automorphisms, then its 
automorphism
group is a so-called {\sl Hurwitz group}, i.e., a finite quotient of the triangle group $\Delta(2,3,7)$. Equivalently, the Riemann surface is an {\'e}tale cover of the Klein quartic $X(7)$. However, it is hitherto unknown which finite groups
can occur as Hurwitz groups. It is even true that, for every integer 
$n$, there exists a $g$ such that there are more than $n$ non-isomorphic
Riemann surfaces of genus $g$ which attain Hurwitz's bound (Cohen, \cite{Cohen:81}). 
In this note, we want to show that the corresponding questions for 
Mumford curves of genus $g$ over non-archimedean valued fields of positive 
characteristic have a very easy answer: the maximal automorphism groups
can be explicitly described, and they occur for an explicitly given 
1-parameter family of curves (at least for $g \notin \{5,6,7,8\}$). 

The set-up for our result is as follows. 
Let $(k,|\cdot|)$ be a non-archimedean valued field of positive characteristic, and $X$ a Mumford curve (\cite{Mumford:72}, \cite{Gerritzen:80})
of genus $g$ over $k$. This means that
the stable reduction of $X$ over the residue field $\bar{k}$ of $k$ is a union of rational curves intersecting 
in $\bar{k}$-rational points. Equivalently, as a rigid analytic space over $k$, the analytification
$X^{\an}$ of $X$ is isomorphic to an analytic space of the form $\Gamma \backslash ({\bf P}^{1,\an}_k-\mathcal{L})$,
where $\Gamma$ (the so-called Schottky group of $X$) is a discrete free subgroup of $PGL(2,k)$ of rank 
$g$ (acting in the obvious way on ${\bf P}^{1,\an}_k$) with $\mathcal{L}$ as its set of limit points.

In \cite{Cornelissen:01}, it was shown that a Mumford curve of genus $g\geq 2$ satisfies
\begin{equation} \label{bound} |\Aut(X)| \leq \max \{ 12(g-1), F(g) \}, \mbox{ where }  F(g):=2 \sqrt{g} (\sqrt{g}+1)^2.
\end{equation}

\begin{remark}
Note that $F(g) \leq 12(g-1)$ precisely when $g \in \{5,6,7,8\}$. If 
the genus is in this range, it seems to be a difficult (though
manageable) task to find all curves 
with $12(g-1)$ automorphisms, but the problem is of a different nature
than the one under consideration here. Let us mention that there
are at least three one-dimensional families of curves of genus $6$ 
that have $12(g-1)$ automorphisms (\cite{Kato:01}).
\end{remark}

For $g \notin \{5,6,7,8\}$, it was shown in \emph{loc.\ cit.} that the so-called \emph{Artin-Schreier-Mumford} curves $X_{t,c}$
attain the bound (\ref{bound}). Here, $t$ is a non-negative integer, $c \in k^*$ satisfies $|c|<1$, and the 
affine equation of $X_{t,c}$ is given by 
\begin{equation} X_{t,c} \ : \ (y^{p^t}-y)(x^{p^t}-x)=c. \end{equation}
The genus of $X_{t,c}$ is $g_t:=(p^t-1)^2$, and its automorphism group is isomorphic to 
$A_t:={\bf Z}_p^{2t} \sd D_{p^t-1}$, where ${\bf Z}_n$ denotes the cyclic group of order $n$, and $D_n$
is the dihedral group of order $2n$. The action of the automorphisms on the coordinates is
explicitly given by interchanging $x$ and $y$, adding an element of ${\bf F}_{p^t}$ to $x$ or $y$,
and multiplying $x$ and $y$ with $a \in {\bf F}^*_{p^t}$ and $a^{-1}$ respectively.
For fixed $t$  and varying $c$, the family $X_{t,c}$ over the punctured unit disc
$\Delta^*:=\{0<|c|<1\}$ is non-constant (since it extends to a semistable family over the unit disk with a singular fiber at $c=0$). Actually, 
the curve $X_{t,c}$ is exactly Mumford if and only if $|c|<1$, since otherwise, its reduction 
is irreducible over $\bar{k}$. The aim of this
paper is to show the following:

\begin{theorem}
Let $g \notin \{5,6,7,8\}$. Every Mumford curve over a non-archimedean valued field $k$ of positive characteristic
with the maximal number $F(g)$ of automorphisms is isomorphic to an Artin-Schreier-Mumford curve, i.e., 
an element of the family $\{X_{t,c}\}_{c \in \Delta^*}$ for $t=\log_p (\sqrt{g}+1)$, which, in particular, has to
be an integer. 
\end{theorem}

\section{Proof of the theorem}

Let $X$ be a Mumford curve of genus $g$ with $F(g)$ automorphisms. Remark that $\Aut(X)=N/\Gamma$,
where $\Gamma$ is the Schottky group of $X$ and $N$ is its normalizer in $PGL(2,k)$ (\cite{Gerritzen:80}, VII.1). 
Let us denote $A=\Aut(X)$, and recall that the only information given about $A$ is that its
order is $2\sqrt{g}(\sqrt{g}+1)^2$, where $g$ is the rank of $\Gamma$. Let $\bar{\cdot} \ : \ N \rightarrow A$ 
denote the reduction map modulo $\Gamma$.

For integers $t,n$ such that $n|p^t-1$, let $B(t,n)={\bf Z}_p^t \sd {\bf Z}_n$, where the semi-direct action is
that of the Borel subgroup of $PGL(2,k)$. Let $D_n = {\bf Z}_n \sd {\bf Z}_2$ denote the dihedral group of order $2n$. 

The first part of the proof consists in showing the following, which is 
the (much easier) non-archimedean analogue of finding all Hurwitz groups:

\begin{proposition} \label{prop1}
If $X$ is a Mumford curve of genus $g \notin \{5,6,7,8\}$ with $F(g)$ automorphisms, then its genus is of the form 
$g=(p^t-1)$ for some integer $t$, its automorphism group
$A$ is isomorphic to ${\bf Z}_p^{2t} \sd D_{p^t-1}$,  $A\backslash X={\bf P}^1$
and  $X \rightarrow A \backslash X$ is ramified above 3 points with ramification groups $({\bf Z}_2, {\bf Z}_2, B(t,n))$ if $p \neq 2$ and ramified above
2 points with groups $({\bf Z}_2, B(t,n))$ if $p=2$.
\end{proposition}

\begin{proof} Let $\ast$ denote amalgamation product.
The computations in paragraph 6 of \cite{Cornelissen:01} (cf.\ (6.9)--(6.12)) show the following:

\begin{lemma}
If $|A|=F(g)$ then, as abstract groups, $N \cong B(t,n) \ast_{{\bf Z}_n} D_n$ for $n=p^t-1$, and
$g=(p^t-1)^2$.
Furthermore, if $\gamma$ denotes an involution in $D_n - {\bf Z}_n$,
$E$ denotes ${\bf Z}_p^t$ \textup{(}seen as a normal subgroup of $B(t,n)$\textup{)}, and $E'=\gamma E \gamma$
then $\bar{E} \cap \bar{E}' = \{ 1 \}$. 
 
\end{lemma}

Now observe that, since $\Gamma$ is free, any finite subgroup of $N$ is mapped by $\bar{\cdot}$ to an isomorphic
image. Also, the image of any relation in $N$ holds in $A$ -- this applies in particular
to the relations given by the semi-direct product structure in $B(t,n)$ and $D_n$.
In the end, $A$ is a group of order $2 p^{2t}(p^t-1)$
which contains the following subgroups: two elementary abelian groups $\bar{E}, \bar{E}'$ which 
are disjoint, both of which are normalized by a cyclic group $\bar{{\bf Z}}_n$ of order $p^t-1$, which
in its turn is normalized by an element $\bar{\gamma}$ of order 2. 
Now $A$ is generated by the {\sl disjoint} groups $\langle \bar{E}, \bar{E}', \bar{{\bf Z}}_n,
\bar{{\bf Z}}_2 \rangle $, and the order of $A$ equals the product of 
the order of these groups. Hence $\bar{E}$ and $\bar{E}'$ generate
a subgroup of $A$ (namely, $\bar{E} \cdot \bar{E}'$) of order $p^{2t}$.
We will now rely on the
following group-theoretical lemma:

\begin{lemma}
Let $q=p^t$ be a power of a prime number, $G$ a group of order $q^2$, generated by two disjoint
subgroups $\mathcal{E},\mathcal{E}'$ of order $q$. Assume that ${\bf Z}_{q-1}$ acts by automorphisms on $G$, stabilizing $\mathcal{E}$ and $\mathcal{E}'$ such that the restriction
of this action to $\mathcal{E}-\{1\}$ and $\mathcal{E}'-\{1\}$ is simply transitive. Then $G=\mathcal{E} \times \mathcal{E}'$.
\end{lemma}

\begin{proof}
Since $G$ is a $p$-group, it has a non-trivial center $Z$. Let 
$z \in Z$ be a non-trivial element. We claim that no element of ${\bf Z}_{q-1}$ fixes $z$ (we will write the action exponentially). 

Indeed, write $z=ab$ for $a \in \mathcal{E}$ and $b \in \mathcal{E}'$. This way
of writing is unique, since $G=\mathcal{E} \cdot \mathcal{E}'$ and $\mathcal{E} \cap \mathcal{E}' = \{1\}$. 
Let $\sigma \in {\bf Z}_{q-1}$, then if $z^\sigma = z$, $ab=a^\sigma b^\sigma $, so $a=a^\sigma $
and $b = b^\sigma$
by the uniqueness of writing. However, since the cyclic group is assumed to act simply transitively
on $\mathcal{E}-\{1\}$ and $\mathcal{E}'-\{1\}$, this can only happen if $a=b=1$, so $z=1$, a contradiction. 

Hence the set the set ${\bf Z}_{q-1} \cdot z $ has $q-1$ elements, and it belongs to $Z$, since ${\bf Z}_{q-1}$ acts by
automorphisms on $G$.

Since $Z$ is a subgroup of $G$, its order has to divide $q^2$, so it is at least $q$. Assume
that $Z \cap \mathcal{E} =\{1\}$. Then $G = Z \cdot \mathcal{E}$, so $|Z|=q$, and since $Z$ is normal in $G$, 
$G=Z \sd \mathcal{E}$, so actually $G=Z \times \mathcal{E}$, contradicting the fact that $Z$ is the center of $G$.

Hence we can find a non-trivial element $\epsilon \in Z \cap \mathcal{E}$, so that $[\epsilon,\mathcal{E}']=1$.
Acting on this with ${\bf Z}_{q-1}$ (which is transitive on $\mathcal{E}$) shows that $[\mathcal{E},\mathcal{E}']=1$, so
that $G=\mathcal{E} \times \mathcal{E}'$ as desired. 
\end{proof}

This lemma implies immediately that $A = (\bar{E} \times \bar{E'}) \sd \bar{D}_n$. Since the 
statement about $A\backslash X$ is in \cite{Cornelissen:01} (cf.\ (2.4)), the proof
of proposition \ref{prop1} is finished. 
\end{proof} 

\begin{remark}
>From the above proof, we see that $\Gamma = \ker(\bar{\cdot})$ contains $[E,\gamma E \gamma]$, and actually 
has to equal it, since $(E \ast E') /[E,E'] = E \times E'$. Up to conjugation, one 
can embed $E$ into $PGL(2,k)$ only as upper triangular matrices 
$$ E = \{ \left( \begin{array}{cc} 1 & V \\ 0 & 1 \end{array} \right) \} $$
for some ${\bf F}_p$-vector space $V$ of dimension $t$ in $k$, and actually, since a cyclic group
${\bf Z}_n$ ($n=p^t-1$) has 
to act semidirectly on it, this cyclic ${\bf Z}_n$ has to be embedded as diagonal matrices over ${\bf F}^*_{p^t}$, and
hence $V$ has to be of the form ${\bf F}_{p^t} \cdot x$ for some $x \in k^*$, as a little matrix
calculation shows. We can assume $x=1$ by conjugating again with a suitable diagonal matrix. 
Then $\gamma$ can only be embedded as 
$$ \gamma = \left( \begin{array}{cc} 0 & 1 \\ C & 0 \end{array} \right) $$
for some $C \in k^*$. Thus, up to $PGL(2,k)$-conjugation (viz., isomorphism of the 
corresponding Mumford curves, cf. \cite{Gerritzen:80}, IV.3.10), the group $\Gamma$ is 
completely characterized by specifying the number $C \in k^*$, which should satisfy $|C|>1$
so that $\Gamma$ is indeed a discrete subgroup of $PGL(2,k)$. The latter fact follows,
e.g., via the method of isometric circles, cf \cite{Cornelissen:01}, \S 8: the isometric
circle of any element of $E$ is (in our normalization) the unit circle $\{|z|=1\}$, whereas
the isometric circle of any element of $E'$ is $\gamma \cdot \{|z|=1\} = \{|z-C|=1\}$. For $E \ast
E'$ to be discretely embedded in $PGL(2,k)$, these isometric circles should not 
intersect, leading to $|C|>1$. 

It would be interesting to find the relation between the parameters $c$
of the algebraic description $X_{t,c}$ 
and $C$ in this description of the Schottky group. 
\end{remark}

The theorem will now follow from the following algebraic fact:

\begin{proposition}
If $X$ is an algebraic curve over a field $k$ of characteristic $p>0$ whose automorphism
group is isomorphic to $A={\bf Z}_p^{2t} \sd D_{p^t-1}$. Assume
that, for $p \neq 2$,  the quotient by $A$ is of the form 
$X \rightarrow A\backslash X={\bf P}^1$
and is ramified above 3 points, say, $(P_1,P_2,P_3)$ with ramification groups $({\bf Z}_2,{\bf Z}_2,B(t,n))$. If $p=2$ suppose that two points $(P_1,P_2)$
are ramified with groups $({\bf Z}_2,B(t,n))$.
Then
$X$ is isomorphic to an Artin-Schreier curve $X_{t,c}$ for some $c \in k^*$.

\end{proposition}

\begin{proof}
Let $p \neq 2$. The automorphism cover $X \rightarrow A \backslash X = {\bf P}^1$ has to decompose into 
a tower of successive Galois extensions as follows:
$$ X \rightarrow X_1:={\bf Z}_p^{2t} \backslash X \rightarrow 
X_2:= {\bf Z}_n \backslash X_1 \rightarrow X_3={\bf Z}_2 \backslash X_2.$$
For this tower to have the correct ramification behaviour, the following should hold.
In $X_2 \rightarrow X_3$, exactly $P_1$ and $P_2$ should ramify, hence $X_2={\bf P}^1$;
let $P_{3,1}$ and $P_{3,2}$ denote the points of $X_2$ above $P_3$. Since $X_1 \rightarrow X_2$
is a tame cover in which only the two points $P_{3,i}$ should ramify, they have to ramify 
completely, and $X_1={\bf P}^1$. Denote the points in $X$ above $P_{3,i}$ in $X_1$ by the same symbol. Now
$X \rightarrow X_1={\bf P}^1$ is an elementary abelian $p$-cover in which two points should 
ramify with ramification
groups ${\bf Z}_p^t$. We can reorder the situation so that $X$ admits two quotients $X^{(1)}=
E_1 \backslash X$
and $X^{(2)}=E_2 \backslash X$ for two groups $E_1=E_2={\bf Z}_p^t$ such that $P_{3,i}$ branches
completely in $X^{(i)}$ and $P_{3,j}$ is unbranched in $X^{(i)}$ for $i \neq j$. Let $z$ be a coordinate
on $X_1={\bf P}^1$ and set $z(P_{3,1})=0$ and $z(P_{3,2})=\infty$. The equation of $X^{(1)}$ 
has to be a succession of $p$-covers of ${\bf P}^1$ in which exactly one point ramifies completely.
If we decompose it into $t$ totally ramified $p$-covers $X^{(1)}_i \rightarrow X^{(1)}_j$,
one can see inductively that each of these is given by an equation
$y_i^p-y_i= y_{i-1}$ with $y_1=z$ (by suitable normalization), which implies inductively that all $X^{(1)}_i$
are isomorphic to ${\bf P}^1$, and in the end 
we find that $X^{(1)}$ is given by $x^{p^t}-x=z$. A similar argument applies to $X^{(2)}$, 
but an equation of the form $y^{p^t}-y=\frac{c}{z}$ comes out for some $c \in k^*$ (since $\infty$
is supposed to ramify in $X^{(2)}$). Note that we cannot normalize $c=1$, since the
coordinate $z$ was already normalized by the tower of $X^{(1)}$.  Finally, the fiber product 
\begin{equation} \label{fiber}
X^{(1)} \times_{X_1} X^{(2)} 
\end{equation}
dominates $X$, and it has the same degree over $X_1$ as $X$, hence it equals $X$. But the
curve (\ref{fiber}) is exactly equal to $X_{t,c}$, which finishes the proof for
$p \neq 2$.  

For $p = 2$, except for the fact that $X_2 \rightarrow X_3$ is 
Artin-Schreier of order two with a unique totally ramified point (instead of Kummer with two ramified points), the proof is entirely analogous, and we refrain from presenting
the details.  
\end{proof}

\begin{remark}
If the moduli space of Mumford curves of genus $g$ over $k$ is stratified according
to automorphism groups, then this proposition shows that the stratum with maximal
automorphism group is exactly equal to the locus of Artin-Schreier-Mumford curves, 
and hence it is rigid analytically connected. This connectedness statement can fail
to hold for more general strata, see for example \cite{Kato:01}. Nevertheless, the dimension
of more general equivariant first order deformation spaces of Mumford curves can be computed in term of the ramification data associated to 
the automorphism group (or the tree of groups associated with the normalizer of the 
Schottky group), cf. \cite{Cornelissen:02}.
\end{remark}

\bibliographystyle{amsplain}

\end{document}